\documentclass[a4paper,10pt,twocolumn]{article}

\usepackage{amsthm,amsmath,amsfonts,amssymb}
\usepackage{graphicx}

\newtheorem{lemma}{Lemma}[section]
\newtheorem{definition}{Definition}[section]
\newtheorem{theorem}{Theorem}[section]
\newtheorem{corollary}{Corollary}[section]
\newtheorem{example}{Example}[section]

\newtheorem{assumption}{Asumption}[section]

\def\L2{{\cal L}_2}

\def\begce{\begin{center}}
\def\endce{\end{center}}
\def\begar{\begin{array}}
\def\endar{\end{array}}
\def\begeq{\begin{equation}}
\def\endeq{\end{equation}}
\def\begdi{\begin{displaymath}}
\def\enddi{\end{displaymath}}
\def\begdis{\begin{eqnarray*}}
\def\enddis{\end{eqnarray*}}
\def\begeqa{\begin{eqnarray}}
\def\endeqa{\end{eqnarray}}
\def\begdes{\begin{description}}
\def\enddes{\end{description}}
\def\begit{\begin{itemize}}
\def\endit{\end{itemize}}
\def\begen{\begin{enumerate}}
\def\enden{\end{enumerate}}
\def\beglar{\left[\begin{array}}
\def\endrar{\end{array}\right]}
\def\begle{\begin{lemma}}
\def\endle{\end{lemma}}
\def\begde{\begin{definition}}
\def\endde{\end{definition}}
\def\begth{\begin{theorem}}
\def\endth{\end{theorem}}
\def\begco{\begin{corollary}}
\def\endco{\end{corollary}}
\def\begex{\begin{example}}
\def\endex{\end{example}}

\def\begas{\begin{assumption}}
\def\endas{\end{assumption}}
\def\begpp{\begin{property}}
\def\endpp{\end{property}}

\def\L2{{\mathcal L}_2}

\def\R{\mathbb{R}}

\newcommand{\rv}[1]{\boldsymbol{#1}} 
\def\1{\rv{1}}

\def\gti{\rightarrow \infty}
\def\gt{\rightarrow}

\def\L{\emph{\tiny L}}

\title{Towards WaterLab: A Test Facility for New Cyber-Physical Technologies in Water Distribution Networks
}

\usepackage{authblk}
\author[1]{Arturo Tejada Ruiz\thanks{Corresponding author. Email: \{arturotejadaruiz, klaudiahorvath, humbertosteinshiromoto, heddebosman\} @ incas3.eu}}
\author[1]{Klaudia Horv\'{a}th}
\author[1]{Humberto Stein Shiromoto}
\author[1]{Hedde Bosman}
       \affil[1]{INCAS$^3$, Dr. Nassaulaan 9, 9401HJ, Assen, The Netherlands }

\date{15 April, 2015}

\begin{document}

\maketitle
\begin{abstract}
This paper reports the initial steps in the development of WaterLab, an ambitious experimental facility for the testing of new cyber-physical technologies in drinking water distribution networks (DWDN). WaterLab's initial focus is on wireless control networks and on data-based, distributed anomaly detection over wireless sensor networks. The former can be used to control the hydraulic properties of a DWDN, while the latter can be used for in-situ detection and isolation of contamination and hydraulic faults.
\end{abstract}


\section{Introduction}

Drinking water production and distribution are tightly regulated activities. In most countries, water utilities must guarantee minimum levels of water quality by keeping the concentration of particular water contaminants (e.g., bacteria \cite{Blokker.2010}) below levels harmful to consumers \cite{Tagena-et.al.2011}. Water quality is traditionally monitored by analyzing water samples in laboratories every few days or weeks. Although there is growing interest in supplementing this relatively slow process with data acquired from new, real-time water quality sensor systems \cite{paperEliadesPolycarpou,vanThienen.2014}, the testing of this and other new cyber-physical technologies (e.g., micro-turbines for energy recovery \cite{McNabola-et.al.2014} or data-driven anomaly detection \cite{Bosman-et.al.2014}) in working drinking water distribution networks (DWDNs) is severely limited due to safety regulations.

This paper reports the initial steps in the development of WaterLab, an experimental facility specifically designed to address this limitation. WaterLab will be equipped with a flexible, reconfigurable hydraulic network that, although smaller than real DWDNs, will provide the same hydraulic conditions (in terms of time constants, flow patterns, etc.) found in practice. This will be achieved via an advanced hydraulics control system deployed over a state-of-the-art wireless sensor/actuator network composed of hydraulic sensors, motorized valves, and ad-hoc wireless communication nodes. WaterLab is expected to help testing several new cyber-physical systems technologies starting with wireless control networks (WCNs), a control paradigm new to the water community, and data-based, distributed anomaly detection algorithms executed by wireless communication nodes. The facility's main elements and our initial analysis are described in Section 2 next. Section 3 contains our conclusions.

\section{WaterLab's Main Elements}

\subsection{Hydraulic Network}

\begin{figure}[tb] \centering
\includegraphics[width=0.45\textwidth]{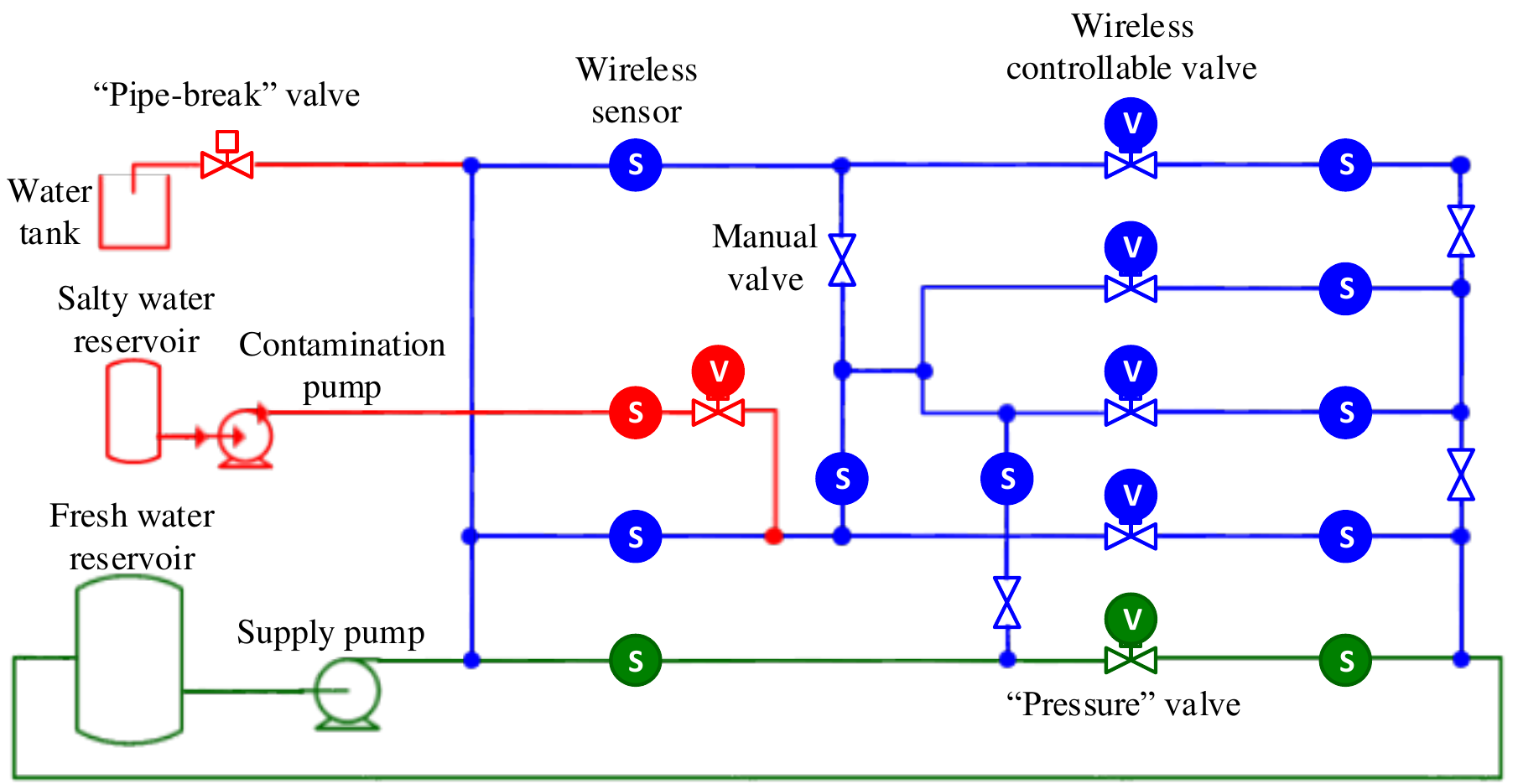}
\caption{\label{fig:schematic} WaterLab's  main features.}
\end{figure}

As the simplified diagram in Figure \ref{fig:schematic} shows, WaterLab's hydraulic network is composed of a supply loop (green, bottom), a DWDN (top, blue), and a set of experimental hydraulic devices (red, left). The supply loop is designed to provide an almost constant water pressure differential to the DWDN by establishing a high water flow rate across a ``pressure'' valve. The DWDN mimics a residential water network. The water flow in each of its branches will be independently adjusted to match the time-varying flow patterns found in typical residential DWDNs. In addition, the DWDN will be equipped with manual valves to allow testing under different DWDN topologies. Finally, the experimental hydraulic devices are add-ons that will allow for, among others, the testing of anomaly detection algorithms for pipe-breaks and contamination faults. The whole setup is expected to be at least 20 meters long and to be stacked vertically over ground.

\begin{figure}
\centering
\includegraphics[width=1\columnwidth]{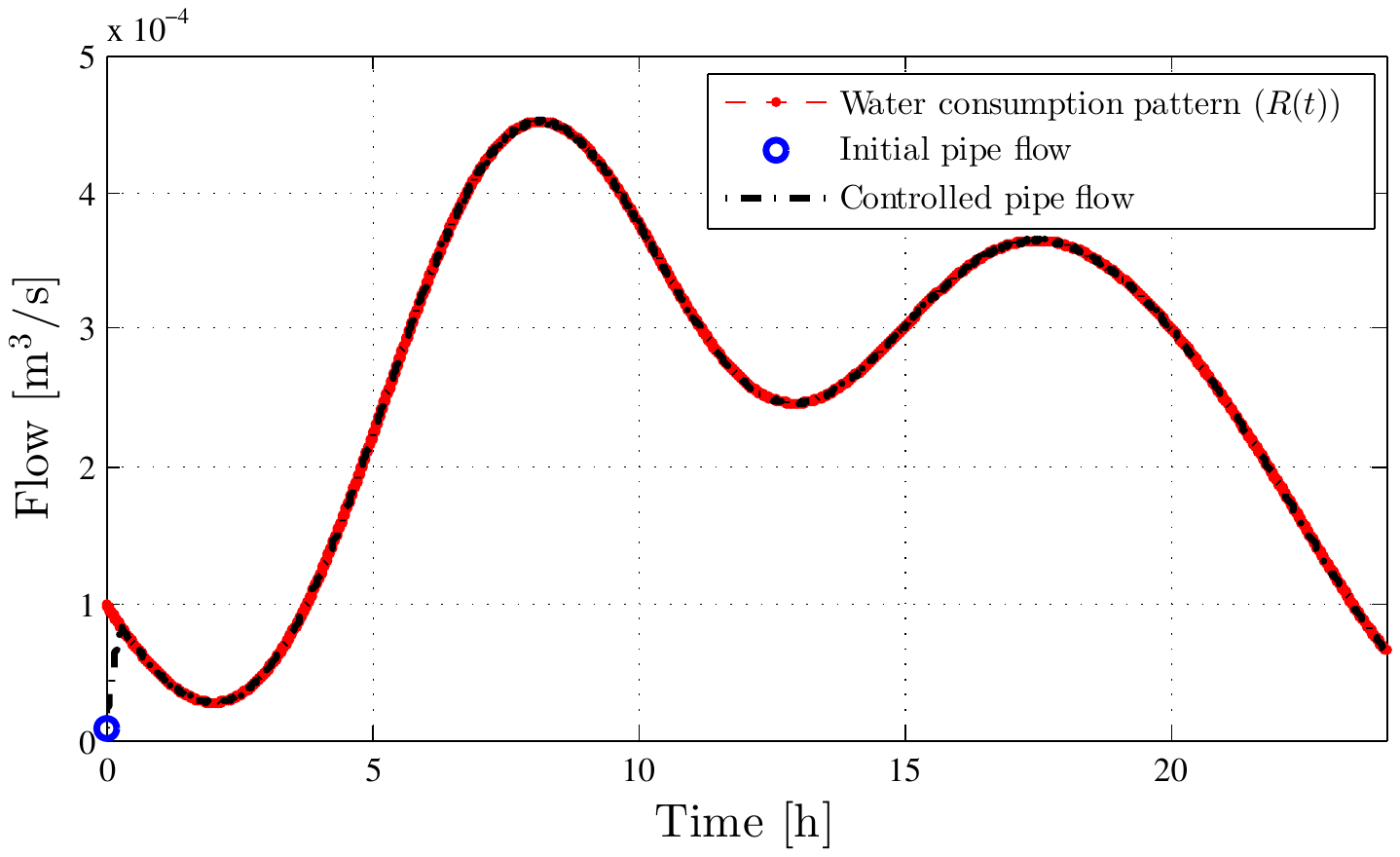}
\caption{Sketch of a typical daily water consumption pattern (red dashed line). The benefits of the proposed water flow control law are also shown (black dot-dashed line).}
\label{fig:DistortedCamel}
\end{figure}

An initial analysis of water consumption data from a water utility's DWDN revealed a daily water consumption pattern similar to that sketched in Figure~\ref{fig:DistortedCamel}. As this figure shows, the highest water consumption occurs during the morning and the evening (when people shower, prepare breakfast, prepare dinner, etc.), and the lowest consumption occurs at night (when most people are asleep). Clearly, there is considerable daily fluctuation. The data also revealed that near the consumption points the most frequent pipe diameters are between 40mm and 50mm, and that the flow presents slow time constants and is mainly laminar. The setup design and modeling (see below) are based on these conclusions.

\subsubsection{Flow Modelling}

Due to the afore-described flow conditions, it was concluded that the water inertial effects, the fluid compressibility and the pipe's elastic behaviour could be disregarded. As a consequence, the flow in the pipes can be accurately modeled by one-dimensional rigid water-column equations \cite{Larock-et.al.2000}, which are derived from the Navier-Stokes equations after certain simplifications (see \cite{Piller2007}) and are extensively used by researchers and practitioners \cite{Epanet2000,Piller2006,dePersis-Kallesoe.2011}. According to this approach, the flow in every pipe equipped with a valve (see Figure \ref{fig:schematic}) can be modeled by
\begin{equation*}\label{eq:one momentum equation general}
\dot q=gA\frac{H_{u}-H_{d}}{L}-\frac{f}{2 D A}q |q| -k_{v}\frac{q |q|}{2 A L},
\end{equation*}
where $q$~(m$^3$/s) denotes the pipe flow, $g$ the acceleration of gravity (m/s$^2$), $A$~(m$^2$) is the cross section of the pipe, $H_{u}$~(m) is the upstream pressure head, $H_{d}$~(m) is the downstream pressure head, $L$~(m) is the pipe's length , $f$ is the friction factor, $D$~(m) denotes the pipe's diameter,  and $k_{v}$ is the valve's head loss coefficient. In this equation, the first term in the right hand side (RHS) denotes the change in flow due to the pressure gradient (assumed to be constant as explained before). The second RHS term is the change in flow due to friction (which is generally a nonliner function of $q$). The last RHS term denotes the pressure loss due to the valve. In this term, $k_{v}$ is a time-variable coefficient that depends on the opening position of the (controllable) valve. Thus, $k_v$ is the control action variable. In the laminar regime $f\propto 1/q$, so the equation above can be simplified as follows:
\begin{equation*}\label{eq:one momentum equation laminar}
\dot q=gA\frac{H_{u}-H_{d}}{L}-\frac{64 \nu }{2 D ^2}q -k_{v}\frac{q |q|}{2 A L},
\end{equation*}
where  $\nu$ is the kinematic viscosity ($m^2/s$). Finally, if the flow is unidirectional (a valid assumption when there is only one source of water pressure, as is the case in WaterLab), the rigid water-column equations further simplify to

\begin{equation}\label{eq:RWC model}
\dot q={\cal K}-{\cal L}q-q^2u,
\end{equation}
where ${\cal K}=gA\frac{H_{u}-H_{d}}{L}>0$, ${\cal L}=\frac{64 \nu }{2 D ^2}>0$ and $u=\frac{k_{v}}{2 A L}$. This equation is well suited for the derivation of hydraulic control laws as explained below (see also \cite{dePersis-Kallesoe.2011}).

\subsection{Wireless Sensor/Actuator Networks}

WaterLab will be equipped with (among others) flow rate, pressure, temperature, conductivity, and turbidity sensors (``S'' in Figure \ref{fig:schematic}) and motorized valves (``V'' in Figure \ref{fig:schematic}). All sensors and valves will be interconnected through a wireless communication network based on an ad-hoc platform (e.g., Telosb \cite{Willow.2014}) that can establish two-way communication using a variety of communication protocols (e.g., ZigBee, gossiping, etc.). The network could then be used to transport sensor data to a centralized location as is done in traditional DWDN decision support systems. Alternatively, the wireless nodes themselves could be used to process the data and implement decentralized algorithms for automatic control or anomaly detection purposes as explained next.

\subsubsection{Wireless Control Network}

A wireless control network (WCN) is a distributed control algorithm executed by the communication nodes of a wireless sensor/actuator network \cite{Pajic-et.al.2011}. WaterLab will use a WCN to establish desired water flow patterns in its DWDN (see Section 2.1). The WCN design requires two steps: First, a control law is derived based on the rigid water column equations using systems and control theory. Second, the control law is  parallelized. That is, its computation is divided into smaller parts that are then distributed among the wireless nodes. Our initial analysis shows that good performance can be obtained by using Sontag's Formula \cite{Lin-Sontag.1991}. This control technique can take into account the valves' limitations (e.g., saturation), attenuates disturbances, and can be easily distributed among the nodes in the form of, for instance, a neural network. A summary of this approach is given next.

\subsubsection{Flow Control via Sontag's formula}

Consider again the typical daily water consumption pattern sketched in Figure \ref{fig:DistortedCamel}. Since it is periodic, it can be well approximated by a truncated Fourier series of the form:
\begin{equation}\label{eq:reference signal}
	R(t)=\sum_{i=1}^3b_i\cos(i\omega t)+a_i\sin(i\omega t),\quad t\geq0,
\end{equation}
where the constants $\omega>0$, and $a_i,b_i\in\mathbb{R}$, $i=1,2,3$, are such that $R(t)>0$ for every $t\geq0$. The purpose of the control system is to force the flow in the pipes to follow the daily flow pattern shown in Figure \ref{fig:DistortedCamel}, by changing the valve positions appropriately.  This idea is made more precise next.\\

\noindent{\bfseries Control Law}\\

Let $e=q-R$ denote the error between the actual and the desired flow in a pipe. It follows from \eqref{eq:RWC model} that the error dynamics are given by
\begin{eqnarray}
		\dot{e}&=&\dot{q}-\dot{R}\nonumber\\
		&=&{\cal K}-{\cal L}q-q^2u-\dot{R}\nonumber\\
		&=&{\cal K}-{\cal L}(e+R)-\dot{R}-(e+R)^2u\nonumber\\
		&=&{\cal K}-{\cal L}(e+R)-\dot{R}-(e+R)^2\bar{u}-(e+R)^2v,\nonumber\\
		&:=&f(t,e) + g(t,e)v\label{eq:error system}
\end{eqnarray}
where $v=u-\bar{u}$, for some average valve opening position $\bar{u}$. The control law $v$ is a real function of time, $t$, and of $e$  (i.e., $v:\R^+\times \R\gt \R$) and should be such that $e(t)\gt 0$ as $t\gti$ for all initial error values $e(0)$. Controls laws that satisfy the aforementioned conditions are called \emph{continuous stabilizing control laws} and are (in general) difficult to implement in practice. A less conservative approach, is to design a control law such that $|e(t)|\gt \epsilon$ as $t\gti$ for an arbitrarily small $\epsilon>0$, provided that $|e(0)|<\delta$, where $\delta$ can be chosen to be arbitrarily large. That is, provided that the initial error is bounded, such \emph{continuous practical stabilizing control laws} can drive the error arbitrarily close to zero over time (see \cite{Isidori:1999}).\\

\noindent{\bfseries Sontag's Formula}\\

Consider the control Lyapunov function \cite{Jiangetal:2009}
\begin{align*}
V:\mathbb{R}^+\times\mathbb{R}&\to\mathbb{R}^+\\
(t,e)&\mapsto e^2/2,
\end{align*}
and let $L_fV\triangleq\tfrac{\partial V}{\partial e}f$ and $L_gV\triangleq\tfrac{\partial V}{\partial e} g(t,e)$. It can be shown that
\begin{align}\label{eq:Sontag's TV formula}
		v=\begin{cases}
			-\dfrac{L_fV+\sqrt{L_fV^2 + L_gV^4}}{L_gV}&, L_gV(t,e)\neq0\\
			\quad\quad \quad\quad\quad\quad0&,\text{otherwise.}
		\end{cases}
\end{align}
is a continuous stabilizing feedback control law for system \eqref{eq:error system} \cite{Lin-Sontag.1991}. Since this control law is continuous, the closed-loop system can attenuate disturbances such as small measurement noise, small actuator errors, small external disturbances, and small modeling errors (see \cite[Lemma I.2]{SontagWang1996}, also  \cite{FreemanKokotovic2008,Ledyaev:1999} and references therein). In addition, note that due to the Weierstrass Approximation Theorem, \eqref{eq:Sontag's TV formula} can be arbitrarily approximated using (sufficiently high order) polynomials over $e$ and $t$. Moreover, it can be shown that such an approximated controller is a continuous practical stabilizing control law.\\

\noindent{\bfseries Example}\\

Figure \ref{fig:DistortedCamel} shows the result of applying the flow controller to a 100m long pipe under a (constant) head differential of 1.85m, assuming a pipe diameter of 46mm and water viscosity of $10^{-6}$~m$^2$/s. Note that the controller is able to quickly steer the almost zero initial flow ($8.6\times 10^{-6}$ m$^3$/s) to the value required by the reference signal $R(t)$.

\subsubsection{Data-Based, Distributed Anomaly Detection}

The wireless nodes can also be used to detect unexpected behaviors or events, such as leaks or increased contamination, commonly referred to as anomalies \cite{Chandola-et.al.2009}. There is an extensive body of research on this subject, and several anomaly detection methods are already available to recognize unexpected events in sensor data. Anomaly detection is commonly performed at centralized storage and processing facilities. It entails creating models of the sensor data through either mathematical or statistical mean, and then using these models to forecast the behavior of future data points. Discrepancies between the new data and their forecasts indicate the presence of possible anomalies.
Examples of these methods can be found in data-mining systems \cite{phua2010comprehensive,kao2009motivating}, data center management \cite{wang2011statistical}, security card access \cite{biuk2012behavior}, or network intrusion detection \cite{modi2013survey}. 

The use of centralized method has several drawbacks. For instance, in a wireless sensor network, all measured data would need to be transmitted to a selected centralized location. Hence, there is growing interest in distributing anomaly detection methods in embedded network nodes. This approach can provide advantages in terms of both data reduction (e.g., only possible anomalous data points are transmitted) and resource efficiency (e.g., lower bandwidth communication channels). 


WaterLab will provide the opportunity to test several new distributed anomaly detection algorithms based on machine learning (ML) techniques. Particular attention will be given to online recursive regression techniques that are able to fit both linear \cite{bosman2013anomaly} and non-linear \cite{Bosman-et.al.2014} models to sensor data using limited-resource embedded systems. These techniques do not require supervision (i.e., they run autonomously) and are well suited to process data from environments or processes for which little a priori knowledge is available (e.g., the difference between ``normal'' and ``anomalous'' behavior). In such cases, techniques based on learning or self-adaptation are easier to apply than their model-driven counterparts.

Since its operating conditions will be well-known and modelled, WaterLab will provide the opportunity to compare the relative benefits of both data-based and model-based anomaly detection methods. Moreover, WaterLab will also be used to test for the first time in the DWDN context anomaly detection methods (currently under development) based on locally-combined sensor data from multiple sensors nodes.

\section{Conclusions and Next Steps}

This paper summarizes the main features of WaterLab, an ambitious test facility composed of a scaled-down but dynamically equivalent hydraulic network together with a wireless sensor/actuator network (WSAN). WaterLab is designed to carefully mimic the working conditions of real drinking water distribution networks (DWDNs), thus allowing the testing of new cyber-physical technologies that could otherwise not be tested due to safety regulations in commercial DWDNs.

WaterLab initial analysis phase is now completed. The next steps of the project include the careful recreation of WaterLab's hydraulic network (together with its control system) in a real-time simulation environment (e.g., Simulink) and the construction of a prototype hydraulic network (including a supply loop, a simple DWDN, and a small WSAN).

\bibliographystyle{abbrv}
\bibliography{ATRBIB2015}  

\end{document}